\documentclass{amsart}
\usepackage{amsfonts}
\usepackage{amssymb}
\usepackage{amsmath}

\newtheorem{theorem}{Theorem}[section]

\newtheorem{lemma}[theorem]{Lemma}
\newtheorem{proposition}[theorem]{Proposition}
\newtheorem{corollary}[theorem]{Corollary}
\newtheorem{fact}[theorem]{Fact}
\theoremstyle{definition}
\newtheorem{definition}[theorem]{Definition}

\theoremstyle{remark}
\newtheorem{remark}[theorem]{Remark}
\numberwithin{equation}{section}

\input{tcilatex}

\begin{document}
\title{The dynamical Borel-Cantelli lemma and the waiting time problems}
\author{Stefano Galatolo}
\address{Dipartimento di matematica applicata, via Buonarroti 1, Pisa}
\email{galatolo@dm.unipi.it}
\author{Dong Han Kim}
\address{Department of Mathematics, The University of Suwon, Hwaseong
445-743, Korea}
\email{kimdh@suwon.ac.kr}
\thanks{The second author is supported by the Korea Research Foundation
Grant funded by the Korean Government(MOEHRD). (KRF-2005-214-C00178)}

\begin{abstract}
We investigate the connection between the dynamical Borel-Cantelli and
waiting time results. We prove that if a system has the dynamical
Borel-Cantelli property, then the time needed to enter for the first time in
a sequence of small balls scales as the inverse of the measure of the balls.
Conversely if we know the waiting time behavior of a system we can prove
that certain sequences of decreasing balls satisfies the Borel-Cantelli
property. This allows to obtain Borel-Cantelli like results in systems like
axiom A and generic interval exchanges.
\end{abstract}

\maketitle

\today

\section{Introduction}

Let $\{A_{n}\}$ be a sequence of subsets in a probability space $(X,\mu )$.
The classical Borel-Cantelli lemma states that:

\begin{enumerate}
\item if $\sum \mu (A_{n})<\infty $, then $\mu (\lim \sup A_{n})=0,$ that
is, the set of points which are contained in infinitely many $A_{n}$ has
zero measure.

\item Moreover, if the sets $A_{n}$ are independent, then $\sum \mu
(A_{n})=\infty $ implies that $\mu (\lim \sup A_{n})=1$.
\end{enumerate}

Now, let us consider a dynamical system $(X,T,\mu )$ and suppose that $%
T:X\rightarrow X$ preserves $\mu $. In this case, if the sets $T^{-n}A_{n}$
are independent and $\sum \mu (A_{n})=\infty $, then the set of points such
that $T^{n}x\in A_{n}$ infinitely many times as $n$ increases has full
measure.

In a chaotic, mixing measure preserving dynamical system, sets of the form $A
$ and $T^{-n}A$ tend to "behave" as independent in a certain sense, as $%
n\rightarrow \infty $. By this it is reasonable to ask if a statement like
point 2 above is valid. The answer is that it is not always valid. In \cite%
{F} it is shown an example of a mixing system, where the BC property does
not hold even for nice set sequences (decreasing sequences of balls with the
same center). Hence some stronger requirements are needed (some stronger
form of mixing or some stronger constraint in the sequence of sets).

In this context a decreasing sequence $A_{n}$ is said to be a Borel-Cantelli
sequence (BC) if $\sum \mu (A_{n})=\infty $ and $\mu (\lim \sup
T^{-n}A_{n})=1$ (here $\lim \sup S_{n}$ is the set of points which belongs
to infinitely many $S_{n}$). In a mixing dynamical system, thus, the
\textquotedblleft abundance" of BC sequences can be interpreted as an aspect
of strong chaos and stochastic behavior of the system. Indeed is proved (see
e.g. \cite{CK},\cite{Ph},\cite{Dol},\cite{kim2} \cite{FMP}, \cite{T}, \cite%
{G07}) that in many kind of (more or less) hyperbolic or "fast" mixing
systems, various sequences of geometrically nice sets have the BC property.
The kind of sets which are interesting to be considered in this kind of
problems are usually decreasing sequences of balls with the same center
(these are also called shrinking targets, this approach has relations with
the theory of approximation speed, see \cite{HV}, \cite{K}) or cylinders.

Let us consider another concept wich as we will see is closely related to
the Borel-Cantelli property: the waiting time. Let $A$ be a subset of $X$,
let 
\begin{equation*}
\tau _{A}(x)=\min \{n\in \mathbf{N}:T^{n}(x)\in A\}
\end{equation*}%
this is the time needed for $x\in X$ to enter for the first time in $A$. It
is clear that in an ergodic system (when $A$ has positive measure) $\tau
_{A}(x)$ is almost everywhere finite. Intuitively, when $A$ is smaller and
smaller, then $\tau _{A}(x)$ is bigger and bigger. If the behavior of the
system is chaotic enough, one could expect that for most points $\tau
_{A}(x)\sim \frac{1}{\mu (A)}.$ More precisely, let $B(y,r_{n})$ be a
sequence of balls with center $y$ and radius $r_{n}$. We say that $x$
satisfies the waiting time problem (with respect to the sequence of sets $%
B(y,r_{n})$) if 
\begin{equation}
\lim_{n\rightarrow \infty }\frac{\log \tau _{B(y,r_{n})}(x)}{-\log \mu
(B(y,r_{n}))}=1.  \label{sca}
\end{equation}%
In this case if the local dimension$\footnote{%
If $X$ is a metric space and $\mu $ is a measure on $X$ the upper local
dimension at $x\in X$ is defined as $\overline{d}_{\mu }(x)=%
\mathrel{\mathop{\limsup} \limits_{r\rightarrow 0}}\frac{\log (\mu (B(x,r)))%
}{\log (r)}=%
\mathrel{\mathop{\limsup} \limits_{k\in {\bf N}, k\rightarrow
\infty}}\frac{-\log (\mu (B(x,2^{-k})))}{k}$. The lower local dimension $%
\underline{d}_{\mu }(x)$ is defined in an analogous way by replacing $%
\limsup $ with $\liminf $. If $\overline{d}_{\mu }(x)=\underline{d}_{\mu }(x)
$ we denote with $d_{\mu }(x)$ the local dimension at $x.$}$ of $\mu $ at $y$
is $d_{\mu }(y)$, then the measure of balls scales, for small $r$ as $\mu
(B(y,r_{n}))\sim r_{n}^{d_{\mu }(y)}$ and then 
\begin{equation}
\tau _{B(y,r_{n})}(x)\sim r_{n}^{-d_{\mu }(y)}.  \label{dim}
\end{equation}%
A result of this kind has been proved in various kind of chaotic systems (%
\cite{G}, \cite{G05},\cite{G07}, \cite{kimkim} and see also \cite{Non2003}
for cases where it does not hold). Moreover this problem is related to the
distribution of return times (the property holds when the distribution of
return times in small balls tends to be exponential, see \cite{G05}, see
also \cite{L} for other general relations between waiting time and
recurrence time distribution). While in the literature results on
Borel-Cantelli and Waiting time are somewhat similar (and sometime used
togheter, as in \cite{FMP}), as far as we know, no explicit general
relations about these two concepts are stated.

In this note we show that the Borel-Cantelli property and the waiting time
problem are in general strictly connected: in Section~\ref{sec1} we show
that in systems where decreasing sequences of balls have the BC property,
then the waiting time problem is satisfied for almost all points of the
space (Theorem~\ref{bcwt} or in a different point of view, Theorem~\ref{27}%
). In Section~\ref{sec2} we see that there are examples of systems where the
waiting time problem is satisfied but certain decreasing sequences of balls
does not satisfy the BC condition. This says that this kind of
Borel-Cantelli condition is stronger than the one imposed by the waiting
time problem. However, we see in Section~\ref{sec2} that if we impose
further conditions on the sequences of balls we consider, making radius to
decrease in a \textquotedblleft controlled\textquotedblright\ way, then we
have converse statements (Theorem~\ref{mst} and following). This allows to
use results on the waiting time problem to obtain Borel-Cantelli like
results on certain decreasing sequences of balls in systems like axiom A and
generic interval exchanges.

\section{Borel-Cantelli implies Waiting time\label{sec1}}

We assume that $T$ is a measure preserving transformation on a metric space $%
(X,\mu ,d)$. We will prove a general result about the waiting time problem
which generalizes an inequality between waiting time and measure of sets
proved in \cite{G05}, allowing sets which are not necessarily balls. Then we
prove that in system where decreasing sequences of balls have the BC
property the inequality becomes equality, and then in such systems the
scaling behavior of the waiting time is the same as the scaling behavior of
the measure of small balls (\ref{sca}).

\begin{proposition}
\label{ineq}Let $B_{n}$ be a decreasing sequence of measurable subsets in $X$
with $\lim_{n}\mu (B_{n})=0$. Then we have 
\begin{equation*}
\liminf_{n\rightarrow \infty }\frac{\log \tau _{B_{n}}(x)}{-\log \mu (B_{n})}%
\geq 1\quad \text{for a.e. }x.
\end{equation*}
\end{proposition}

\begin{proof}
Choose a subsequence $n_i$ as $n_i = \min \{ n \ge 1 : \, \mu(B_{n}) <
2^{-i} \}.$ If $n_i \le n < n_{i+1}$, then we have $\tau_{B_n} (x) \ge
\tau_{B_{n_i}} (x) $ for every $x$ and $2^{-i} > \mu(B_{n_i}) \ge \mu(B_n)
\ge 2^{-i-1}$. Therefore, if $n_i \le n < n_{i+1}$, for every $x$ 
\begin{equation*}
\frac{\log \tau_{B_{n}}(x)}{-\log \mu (B_{n})} \ge \frac{\log
\tau_{B_{n_i}}(x)}{-\log \mu(B_{n_{i}})} \cdot \frac{-\log \mu(B_{n_{i}})}{%
-\log \mu(B_{n})} > \frac{\log \tau_{B_{n_i}}(x)}{-\log \mu(B_{n_{i}})}
\cdot \frac{ i }{ i+1 },
\end{equation*}%
which implies that 
\begin{equation*}
\liminf_{i\rightarrow \infty }\frac{\log \tau _{B_{n_i}}(x)}{-\log
\mu(B_{n_{i}})} =\liminf_{n\rightarrow \infty }\frac{\log \tau _{B_{n}}(x)}{%
-\log \mu (B_{n})} \quad \text{ for every } x.
\end{equation*}
Thus we may assume that $\mu (B_{n})\leq 2^{-n}$.

Let 
\begin{equation*}
E_{n}=\{x:\frac{\log \tau _{B_{n}}(x)}{-\log \mu (B_{n})}<1-\delta \}
\end{equation*}%
for some $\delta >0$. Then we have 
\begin{equation}  \label{eqsop}
\begin{split}
\mu (E_{n})& =\mu (\{x:\tau _{B_{n}}(x) < \mu (B_{n})^{-(1-\delta )}\}) \\
& =\sum_{1\leq i<\mu (B_{n})^{-(1-\delta )}}\mu (\{x:\tau_{B_{n}}(x)=i\}) \\
& \leq \sum_{1\leq i < \mu (B_{n})^{-(1-\delta )}}\mu (T^{-i}B_{n}) \\
& \leq \mu (B_{n})^{-(1-\delta )}\mu (B_{n})=\mu (B_{n})^{\delta
}<2^{-n\delta }.
\end{split}%
\end{equation}%
Hence $\sum_{n}\mu (E_{n})<\infty $ and by the Borel-Cantelli Lemma we have 
\begin{equation*}
\mu (\limsup E_{n})=0.
\end{equation*}%
Since $\delta $ is arbitrary we have for almost every $x$ 
\begin{equation*}
\liminf_{n}\frac{\log \tau _{B_{n}}(x)}{-\log \mu (B_{n})}\geq 1.
\end{equation*}
\end{proof}

\begin{remark}
In the above proof (in eq. (\ref{eqsop})), 
for each $\epsilon >0$, we can replace $\delta $ with a sequence $\delta
_{n}\rightarrow 0$ such that $\delta _{n}\leq \frac{(1+\epsilon )\log n}{%
-\log \mu (B_{n})}$. The proof is still valid and we have an estimation on
the "speed of convergence" to this limit inequality. Indeed we obtain that
if $B_{n}$ is a decresing sequence and $\mu (B_{n})\leq 2^{-n},$ typical
points will eventually satisfy 
\begin{equation*}
\frac{\log \tau _{B_{n}}(x)}{-\log \mu (B_{n})}\geq 1-\frac{(1+\epsilon
)\log n}{-\log \mu (B_{n})}
\end{equation*}%
for each $\epsilon >0.$ This is interesting when waiting time is used to
give numerical estimations on the local dimension of attractors. Indeed, by
the above result, working like in (\ref{dim}) we see that in general systems
the scaling behavior of the waiting time gives an upper bound to the local
dimension. This can suggest a numerical method to estimate such a dimension
(see \cite{G05}, \cite{CG06}). This remark, hence, gives also an estimation
on the speed this upper bound is approached. This is very general and does
not require assumptions on the system we consider.
\end{remark}

A sequence of sets $A_{n}$ is said to be strongly Borel-Cantelli if in some
sense the preimages $T^{-n}A_{n}$ covers the space uniformly:

\begin{definition}
Let $1_{A}$ be the indicator function of the set $A.$ The sequence of
subsets $A_{n}\subset X$ is said to be a strongly Borel-Cantelli sequence
(SBC) if for $\mu -$a.e. $x\in X$ we have as $N\rightarrow \infty $ 
\begin{equation*}
\frac{\sum_{n=1}^{N}1_{T^{-n}A_{n}}(x)}{\sum_{n=1}^{N}\mu (A_{n})}%
\rightarrow 1.
\end{equation*}
\end{definition}

As mentioned in the introduction, next theorem says that in Borel-Cantelli
systems we have a relation between waiting time and scaling behavior of the
measure of small balls (as in the waiting time problem $\tau_{B(y,r_{n})}(x)%
\sim r_{n}^{-d_{\mu }(y)}$).

\begin{theorem}
Assume that there is no atom in $X$. \label{bcwt} (i) If every decreasing
sequence of balls in $X$ with the same center is BC, then for every $y$ we
have 
\begin{equation*}
\liminf_{r\rightarrow 0}\frac{\log \tau _{B(y,r)}(x)}{-\log \mu (B(y,r))}%
=1\quad \text{for a.e. }x.
\end{equation*}%
(ii) Suppose that $B(y,r)=\{x : d(x,y) \leq r \}$ is the closed ball. If
every decreasing sequence of balls in $X$ with the same center is SBC, then
for every $y$ we have 
\begin{equation*}
\lim_{r\rightarrow 0}\frac{\log \tau _{B(y,r)}(x)}{-\log \mu (B(y,r))}%
=1\quad \text{for a.e. }x.
\end{equation*}
\end{theorem}

\begin{proof}
(i) Fix $y\in X$. Since $y$ is not an atom, we have $\mu (B(y,r))\downarrow
0 $ as $r\rightarrow 0$.

For each positive integer $i$ define $m (i)$ by the smallest positive
integer such that 
\begin{equation*}
\{ r >0 : \frac 1{m(i)+1} \le \mu (B(y,r)) < \frac 1{i} \} \ne \emptyset.
\end{equation*}
Then $m(1) \le m(2) \le \dots $. Choose $r^{\prime}_i$ as 
\begin{equation*}
\frac 1{m(i)+1} \le \mu (B(y,r^{\prime}_i)) < \frac 1{m(i)}, \quad
i=1,2,\dots.
\end{equation*}
Then there is a sequence $i^{\prime}_{k}$ such that $m(i^{\prime}_k) =
i^{\prime}_k$ and $i^{\prime}_{k} \ge 2 i^{\prime}_{k-1}$. Hence we have 
\begin{equation*}
\sum_{i=1}^\infty \mu (B(y,r^{\prime}_{i})) \ge \sum_{i=1}^\infty \frac{1}{%
m(i)+1} \ge \sum_{k=1}^\infty \frac{i^{\prime}_{k}-i^{\prime}_{k-1}}{%
i^{\prime}_{k}+1} \geq \sum_{k=1}^\infty \frac{i^{\prime}_{k}/2}{%
i^{\prime}_{k}+1} =\infty.
\end{equation*}%
By the BC assumption for almost every $x$, $T^{i}x\in B(y,r^{\prime}_{i})$
for infinitely many $i$'s. Therefore, for almost every $x$ we have $\tau
_{B(y,r^{\prime}_{i})}(x)\leq i\leq 1/\mu (B(y,r^{\prime}_{i}))$ infinitely
many $i$'s. Hence for almost every $x$ 
\begin{equation*}
\liminf_{r \to 0} \frac{\log \tau _{B(y,r)}(x)}{-\log \mu (B(y,r))} \leq 1
\end{equation*}%
for infinitely many $i$'s. The other inequality is obtained by Proposition~%
\ref{ineq}.

(ii) Suppose that there is a $y\in X$ such that 
\begin{equation*}
\limsup_{r\rightarrow 0}\frac{\log \tau _{B(y,r)}(x)}{-\log \mu (B(y,r))}>1
\end{equation*}%
for $x$ ranging in a positive measure set. Choose $m(i)$ as in (i). Then
there is a strictly increasing sequence $\{ i_n \}_{n \ge 1}$ such that $%
m(j) = i_n$ for $i_{n-1} < j \le i_n $. Let 
\begin{equation*}
r_n = \inf \{ r >0 : \frac 1{i_n+1} \le \mu (B(y,r)) < \frac 1{i_n} \},
\quad i=1,2,\dots.
\end{equation*}
Note that $\{ r > 0 : \dfrac 1{i_{n+1}} \le \mu (B(y,r)) < \dfrac 1{i_{n}+1}
\} = \emptyset $. Since 
\begin{equation*}
B(y,r_n) = \{ x : d(x,y) \le r_n \} = \bigcap_{m \ge 1} B(y, r_n + \frac 1m),
\end{equation*}
we have 
\begin{equation*}
\mu(B(y,r_n)) = \lim_{m \to \infty} \mu ( B(y, r_n + \frac 1m )) \ge \frac{1%
}{i_n+1}.
\end{equation*}
If $r_{n} \le r < r_{n-1}$, then we have 
\begin{equation*}
\frac 1{i_{n}+1} \le \mu (B(y,r_{n})) \le \mu (B(y,r)) < \frac 1{i_{n}}
\end{equation*}

Let $A_{k}=B(y,r_{1})$ for $k$, $1 \leq k < i_{1}$ and $A_{k}=B(y,r_{n})$
for $k = i_{n}$, $n = 1,2, \dots$. If $i_{n} < k <i_{n+1}$ for some $n\geq 1$%
, let 
\begin{equation*}
A_{k}=%
\begin{cases}
B(y,r_{n}) & \text{ if }k\leq i_{n}\log (i_{n+1}/i_{n}), \\ 
B(y,r_{n+1}) & \text{ if }k>i_{n}\log (i_{n+1}/i_{n}).%
\end{cases}%
\end{equation*}%
Then for $\log (i_{n+1}/i_{n}) < 1$ we have 
\begin{equation}
\sum_{k=i_{n} + 1 }^{i_{n+1}} \mu (A_{k}) = \sum_{k= i_{n}+1}^{i_{n+1}} \mu
(B(y,r_{n+1})) < \frac{i_{n+1} - i_{n}}{i_{n+1}} = 1 - \frac{i_n}{i_{n+1}} <
\log \frac{i_{n+1}}{i_{n}},  \label{ubound1}
\end{equation}
and for $\log (i_{n+1}/i_{n}) \ge 1$ we have 
\begin{equation}
\begin{split}
\sum_{k=i_{n}+1}^{i_{n+1}}\mu (A_{k}) & =\sum_{k=i_{n}+1}^{\lfloor i_{n}\log
(i_{n+1}/i_{n})\rfloor }\mu ( B(y,r_{n})) + \sum_{k=\lfloor i_{n}\log
(i_{n+1}/i_{n})\rfloor +1}^{i_{n+1}}\mu (B(y,r_{n+1})) \\
& <(\lfloor i_{n}\log \frac{i_{n+1}}{i_{n}}\rfloor -i_{n})\frac{1}{i_{n}}%
+(i_{n+1}-\lfloor i_{n}\log \frac{i_{n+1}}{i_{n}}\rfloor )\frac{1}{i_{n+1}}
\\
&\le (1-\frac{i_{n}}{i_{n+1}})\log \frac{i_{n+1}}{i_{n}} <\log \frac{i_{n+1}%
}{i_{n}},
\end{split}
\label{ubound}
\end{equation}%
where $\lfloor t \rfloor$ is the floor of $t$.

On the opposite side, if $\log (i_{n+1} / i_{n}) < 1$, for any $i_{n} <
j\leq i_{n+1}$ we have 
\begin{equation}
\begin{split}
\sum_{k=i_{n}+1}^{j}\mu (A_{k}) & \geq \frac{j - i_{n}}{i_{n+1}+1} = \frac{1
- i_{n}/i_{n+1}}{1 +1/i_{n+1}} - \frac{i_{n+1}-j}{i_{n+1}+1} \\
&> \frac{1-1/e}{1 +1/i_{n+1}} \log \frac{i_{n+1}}{i_n} - \frac{i_{n+1}-j}{%
i_{n+1}+1} \\
&> (1- \frac 1e - \frac 1{i_{n+1}}) \log \frac{j}{i_n} - \frac{i_{n+1}-j}{%
i_{n+1}+1}.
\end{split}
\label{lbound2}
\end{equation}
The second inequality comes from $(1-1/t) > (1-1/e) \log t$ for $0 < t < e$.
In case of $\log (i_{n+1} / i_{n}) \ge 1$, for any $j$ with $i_{n}\log
(i_{n+1} / i_{n}) < j \leq i_{n+1}$ we have 
\begin{equation}
\begin{split}
\sum_{k=i_{n}+1}^{j} \mu (A_{k}) &\geq (\lfloor i_{n}\log \frac{i_{n+1}}{%
i_{n}} \rfloor -i_{n})\frac{1}{i_{n}+1} +(j-\lfloor i_{n}\log \frac{i_{n+1}}{%
i_{n}} \rfloor )\frac{1}{i_{n+1}+1} \\
&\ge (\frac{i_{n}}{i_{n}+1}- \frac{i_{n}}{i_{n+1}+1} ) \log \frac{i_{n+1}}{%
i_{n}} - \frac{i_{n+1}- j}{i_{n+1}+1} \\
&> ( 1 - \frac{1}{i_{n}+1} - \frac{1}{e}) \log \frac{j}{i_{n}} - \frac{%
i_{n+1}- j}{i_{n+1}+1}.
\end{split}
\label{lbound}
\end{equation}
If $\log (i_{n+1} / i_{n}) \ge 1$ and $i_{n} < j \leq i_{n}\log (i_{n+1} /
i_{n})$, then we have 
\begin{equation}
\sum_{k=i_{n}+1}^{j} \mu (A_{k}) \geq \frac{j - i_{n}}{i_{n}+1} = \frac{%
j/i_n -1}{1 + 1/i_n} > \frac{\log (j/i_n) }{1+1/i_n} > ( 1 - \frac{1}{i_{n}}
) \log \frac{j}{i_{n}}.  \label{lbound1}
\end{equation}


Pick $x$ as $\limsup_{r\rightarrow 0}\frac{\log \tau _{B(y,r)}(x)}{-\log \mu
(B(y,r))}>1$. Then for some $\delta >0$ we have infinitely many $n$'s such
that there exists an $r$ with $r_{n}\leq r<r_{n-1}$ satisfying 
\begin{equation*}
\frac{\log \tau _{B(y,r)}(x)}{-\log \mu (B(y,r))}>1+\delta .
\end{equation*}%
As noticed above 
\begin{equation*}
\frac{1}{i_{n}+1}\leq \mu (B(y,r_{n}))\leq \mu (B(y,r))<\frac{1}{i_{n}}\quad 
\text{for }r_{n}\leq r<r_{n-1}.
\end{equation*}%
We have 
\begin{equation*}
1+\delta <\frac{\log \tau _{B(y,r)}(x)}{-\log \mu (B(y,r))}\leq \frac{\log
\tau _{B(y,r_{n})}(x)}{\log i_{n}}.
\end{equation*}%
which implies that $\tau _{B(y,r_{n})}(x)>{i_{n}}^{1+\delta }$. Note that if 
$\tau _{A_{i_{n}}}(x)=\tau _{B(y,r_{n})}(x)>{i_{n}}^{1+\delta }$, since the $%
A_{j}$ are decreasing then $T^{j}x\notin A_{j}$ for $i_{n}\leq j\leq {i_{n}}%
^{1+\delta }$. Thus there are infinitely many $n$'s such that 
\begin{equation}
S_{i_{n}}(x)=S_{\lfloor {i_{n}}^{1+\delta }\rfloor }(x),  \label{gap}
\end{equation}%
where 
\begin{equation*}
S_{N}(x)=\sum_{n=1}^{N}1_{A_{n}}\circ
T^{n}(x)=\sum_{n=1}^{N}1_{T^{-n}A_{n}}(x).
\end{equation*}%
By the SBC assumption we have 
\begin{equation}
\lim_{N}\frac{S_{N}(x)}{\sum_{n=1}^{N}\mu (A_{n})}=1\quad \text{a.e.}.
\label{SBC}
\end{equation}%
But we have by (\ref{ubound1}) and (\ref{ubound}) 
\begin{equation*}
\sum_{k=1}^{i_{n}}\mu (A_{k})<i_{1}\frac{1}{i_{1}}+\sum_{\ell =1}^{n-1}\log 
\frac{i_{\ell +1}}{i_{\ell }}=\log \frac{i_{n}}{i_{1}}+1
\end{equation*}%
and by (\ref{lbound2}), (\ref{lbound}) and (\ref{lbound1}) for $%
i_{m}<\lfloor {i_{n}}^{1+\delta }\rfloor \leq i_{m+1}$ we have 
\begin{equation*}
\begin{split}
\sum_{k=i_{n}+1}^{\lfloor {i_{n}}^{1+\delta }\rfloor }\mu (A_{k})&
>\sum_{\ell =n}^{m-1}(1-\frac{1}{i_{\ell }}-\frac{1}{e})\log \frac{i_{\ell
+1}}{i_{\ell }}+(1-\frac{1}{i_{m}}-\frac{1}{e})\log \frac{\lfloor {i_{n}}%
^{1+\delta }\rfloor }{i_{m}}-\frac{i_{m+1}-\lfloor {i_{n}}^{1+\delta
}\rfloor }{i_{m+1}+1} \\
& >(1-\frac{1}{i_{n}}-\frac{1}{e})\log \frac{\lfloor {i_{n}}^{1+\delta
}\rfloor }{i_{n}}-1>(1-\frac{1}{i_{n}}-\frac{1}{e})(\delta \log i_{n}-\frac{1%
}{\lfloor {i_{n}}^{1+\delta }\rfloor })-1,
\end{split}%
\end{equation*}%
which contraddicts (\ref{gap}) and (\ref{SBC}). By SBC assumption,\ the set
of $x$ contraddicting (\ref{gap}) and (\ref{SBC}) must have zero measure,
hence $\limsup_{r\rightarrow 0}\frac{\log \tau _{B(y,r)}(x)}{-\log \mu
(B(y,r))}\leq 1$ for almost each $x.$ Combining this result with Proposition~%
\ref{ineq} where the opposite inequality is proved, the proof is complete.
\end{proof}

\begin{corollary}
If every centered, decreasing sequence of balls in $X$ is SBC, Then for
every $y$ we have 
\begin{equation*}
\liminf_{r\rightarrow 0}\frac{\log \tau _{B(y,r)}(x)}{-\log r}=\underbar %
d_{\mu }(y),\quad \limsup_{r\rightarrow 0}\frac{\log \tau _{B(y,r)}(x)}{%
-\log r}=\bar{d}_{\mu }(y)\quad \text{for a.e. }x.
\end{equation*}%
If every decreasing sequence of balls in $X$ is BC and $d_{\mu }(y)$ exists,
then 
\begin{equation*}
\liminf_{r\rightarrow 0}\frac{\log \tau _{B(y,r)}(x)}{-\log r}=d_{\mu
}(y)\quad \text{for a.e. }x.
\end{equation*}
\end{corollary}

The waiting time describes the speed the orbit of a certain point $x$
approaches another point $y$. Another way to consider this kind of questions
is to consider the behavior of limits of the form $\liminf_{n\geq 1}n^{\beta
}\cdot d(y,T^{n}(x)).$ Under this approach, in \cite{Bosher}, Boshernitzan
showed the following quantitative recurrence theorem.

\begin{fact}
Let $(X,\Phi ,\mu ,d,T)$ be a metric measure preserving system. Assume that
for some $\alpha >0$, the Hausdorff $\alpha $-measure $H_{\alpha }$ is $%
\sigma $-finite on $X=(X,d)$. Then for almost all $x\in X$ we have 
\begin{equation*}
\liminf_{n\geq 1}n^{\beta }\cdot d(x,T^{n}(x))<\infty ,\text{ with }\beta =%
\frac{1}{\alpha }.
\end{equation*}%
If, moreover, $H_{\alpha }(X)=0$, then for almost all $x\in X$ 
\begin{equation*}
\liminf_{n\geq 1}n^{\beta }\cdot d(x,T^{n}(x))=0.
\end{equation*}
\end{fact}

By the BC properties of the balls we have an analogous quantitative
approximation theorem for the waiting time. In \cite{BGI}, for general
measure preserving transformations, it is proved that for $\mu $-almost all $%
x\in X$ one has 
\begin{equation*}
\liminf_{n\geq 1}n^{\beta }\cdot d(T^{n}x,y)=\infty \text{ with }\beta >%
\frac{1}{\underbar d_{\mu }(y)}.
\end{equation*}
In the case of Borel-Cantelli systems we have

\begin{theorem}
\label{27}Let $(X,\Phi ,\mu ,d,T)$ be a metric measure preserving system. If
every decreasing sequence of balls in $X$ is BC, then for $\mu $-almost all $%
x\in X$ one has 
\begin{equation*}
\liminf_{n\geq 1}n^{\beta }\cdot d(T^{n}x,y)=0\text{ with }\beta <\frac{1}{%
\underbar d_{\mu }(y)}.
\end{equation*}
\end{theorem}

\begin{proof}
Fix a $y \in X$. By the definition of $\underbar d_{\mu }(y)$, if $\beta <%
\frac{1}{\underbar d_{\mu }(y)}$, then there are infinitely many $n_{i}$'s
such that for any $C>0$ 
\begin{equation*}
\mu (B(y,\frac{C}{n_{i}^{\beta }}))\geq \frac{1}{n_{i}}.
\end{equation*}%
Assume that $n_{i}>2n_{i-1}$. Then we have 
\begin{equation*}
\begin{split}
\sum_{n=1}^{\infty }\mu (B(y,\frac{C}{n^{\beta }}))& \geq \sum_{i=1}^{\infty
}(n_{i}-n_{i-1})\mu (B(y,\frac{C}{n_{i}^{\beta }})) \\
& \geq \sum_{i=1}^{\infty }(n_{i}-n_{i-1})\frac{1}{n_{i}} \\
& \geq \sum_{i=1}^{\infty }\frac{1}{2}=\infty .
\end{split}%
\end{equation*}%
The BC condition implies that $T^{n}x\in B(y,\frac{C}{n^{\beta }})$ for
infinitely many $n$'s. Hence, we have for any $C>0$ 
\begin{equation*}
\liminf_{n\geq 1}n^{\beta }\cdot d(T^{n}x,y)\leq C.
\end{equation*}
\end{proof}

\section{Waiting time and Shrinking targets\label{sec2}}

In the previous section we supposed that sequences of nested balls have the
BC or SBC property and we have seen that this implies strong properties
about waiting time behavior. In this section we will see that the two
concepts are not equivalent and in some sense \textquotedblleft waiting time
is weaker than Borel-Cantelli property\textquotedblright . We also will see
in which direction it is possible to weaken the BC property to have some
converse implication. This direction is very natural: indeed we have to
consider sequences of nice sets, as balls with the same center (shrinking
targets) and such that the sequence of radii decreases in a controlled way.
We remark that this kind of general philosophy (weaker mixing assumption,
stronger requirements on the sets) is similar to the one which is present in
the results of \cite{T}. This remark allows to use results on waiting time
to obtain some Borel-Cantelli results in systems like typical Interval
Exchanges or Axiom A systems.

\begin{definition}
We say that that a system $(X,T,\mu )$ has the shrinking target property
(STP)\ if for any $x_{0}\in X$ any sequences of balls centered at $x_{0}$
has the BC property. Moreover we say that a system has the monotone
shrinking target property (MSTP) if any \emph{decreasing }sequences of balls
centered at $x_{0}$ has the BC property.
\end{definition}

In \cite{F} (see also \cite{Kurzweil}) it is proved that no rotations on the 
$d$-dimensional torus have the STP property and moreover, only rotations
having some particular arithmetical property have the MSTP.

More precisely, let us introduce some notation: consider $\alpha \in \mathbf{%
R}^{d},$ and consider the sup norm $|\alpha |=\sup (|\alpha
_{1}|,...,|\alpha _{d}|).$ Moreover, for $\alpha \in \mathbf{R}$ let us
consider the distance to the closest integer $||\alpha ||=\underset{p\in Z}{%
\inf }|\alpha -p|$ and its generalization on $\mathbf{R}^{d}$: $||\alpha
||=\sup_{i}||\alpha _{i}||.$

\begin{definition}
Let $d\geq 1$, the set 
\begin{equation*}
\Omega ^{d}=\{\alpha \in \mathbf{R}^{d}:\exists C>0\ s.t.\ \forall Q\in 
\mathbf{Z}-\{0\},\ ||Q\alpha ||\geq C|Q|^{-\frac{1}{d}}\}
\end{equation*}
is called set of constant type vectors in $\mathbf{R}^{d}.$
\end{definition}

\begin{theorem}
(\cite{F}) Let $\mathbb{T}^{d}$ be the $d$-dimensional torus. Let us
consider the system $(\mathbb{T}^{d},T_{\alpha },\mu )$ where $T_{\alpha }$
is the translation by a vector $\alpha $ and $\mu $ is the Haar measure on $%
\mathbb{T}^{d}.$ Then we have that $\forall \alpha $ $(\mathbb{T}%
^{d},T_{\alpha },\mu )$ does not have the STP, moreover $(\mathbb{T}%
^{d},T_{\alpha },\mu )$ has the MSTP if and only if $\alpha $ is of constant
type.
\end{theorem}

It is known (\cite{Non2003}) that in dimension 1, for almost every $\alpha $
we have 
\begin{equation*}
\lim_{r\rightarrow 0}\frac{\log \tau _{B(y,r)}(x)}{-\log \mu (B(y,r))}%
=1\quad \text{for a.e. }x,
\end{equation*}%
moreover for every $\alpha $ 
\begin{equation*}
\liminf_{r\rightarrow 0} \frac{\log \tau _{B(y,r)}(x)}{-\log \mu (B(y,r))}%
=1\quad \text{for a.e. }x,
\end{equation*}%
Thus the converse of Theorem~\ref{bcwt} (in particular point (i) ) does not
hold, even if we restrict to decreasing families of balls having the same
centers, as in the shrinking targets framework. (See also \cite{Kimnew}.)

One of the reasons why not all translations have the STP is that the radii
of the balls can decrease in any way. Putting some restriction on this
decreasing rate the STP become equivalent to the waiting time problem.

\begin{theorem}
\label{mst}Let $\{B(y,r_{n})\}$ be a decreasing sequence of centered balls
such that 
\begin{equation*}
\underset{n\rightarrow \infty }{\lim \sup }\frac{\log r_{n}}{-\log n}<\frac{1%
}{\underline{d}_{\mu }(y)}.
\end{equation*}%
If 
\begin{equation}
\underset{r\rightarrow 0}{\lim \inf }\frac{\log \tau _{B(y,r)}(x)}{-\log r}=%
\underline{d}_{\mu }(y)  \label{sop}
\end{equation}%
then $x\in \lim \sup T^{-n}(B(y,r_{n})).$
\end{theorem}

We remark that condition \ref{sop} above is implied by the waiting time
problem and if equation \ref{sop} holds for almost each $x$, then 
\begin{equation*}
\mu (\lim \sup T^{-n}(B(y,r_{n})))=1
\end{equation*}%
and then such a $\{B(y,r_{n})\}$ has the BC property.

\begin{proof}
If $\underset{r\rightarrow 0}{\lim \inf }\frac{\log \tau _{B(y,r)}(x)}{-\log
r}={\underline{d}}_{\mu }(y)$, then there is a sequence $\rho _{n}\downarrow
0$ such that for each small $\epsilon >0,$ $x\in \underset{i\leq \rho _{n}^{-%
{\underline{d}}_{\mu }(y)-\epsilon }}{\cup }T^{-i}(B(y,\rho _{n}))$ for each 
$n.$ If $\underset{n\rightarrow \infty }{\lim \sup }\frac{-\log r_{n}}{\log n%
}=\frac{1}{{d}}<\frac{1}{\underline{d}_{\mu }(y)}$ then, when $m$ is big
enough $r_{m}\geq m^{-1/({d}-\epsilon )}.$ Therefore, if $\epsilon $ is
small enough such that $\epsilon \leq (d-\underline{d}_{\mu })/2$, then $%
\rho _{n}\leq r_{\left\lfloor \rho _{n}^{-\underline{d}_{\mu }-\epsilon
}\right\rfloor }$ eventually, with respect to $n$. Hence, since $r_{m}$ is
decreasing, we have 
\begin{equation*}
x\in \underset{i\leq \rho _{n}^{-\underline{d}_{\mu }-\epsilon }}{\cup }%
T^{-i}(B(y,\rho _{n}))\subset \underset{m\leq \rho _{n}^{-\underline{d}_{\mu
}-\epsilon }}{\cup }T^{-m}(B(y,r_{m})).
\end{equation*}%
This is true for infinitely many $n$ and thus $x\in \lim \sup
T^{-n}(B(y,r_{n})).$
\end{proof}

Since in Axiom A systems it holds that \ref{sop} is verified for typical
points (\cite{G}), this implies that in such systems, decreasing sequences
of balls, verifying the assumptions in theorem \ref{mst} have the BC
property. This extend a result of \cite{Dol} (Theorem 7) which requires the
invariant measure to have a smooth density (but has milder requirements on
the hyperbolicity of the system).

We remark that if we have stronger assumptions on the behavior of \ $\tau
_{B(y,r)}(x)$ we can include other kind of sequences $r_{n}$ and generalize
the above theorem to the following

\begin{proposition}
\label{aboveprop}If for some $x,y\in X$ there is a sequence $\rho
_{n}\downarrow 0$ such that $\tau _{B(y,\rho _{n})}(x)<f(\rho _{n})$, with $%
f:\mathbf{R}^{+}\rightarrow \mathbf{R}^{+}$ be invertible and both $f$ and $%
f^{-1}$ are strictly decreasing. Then for each decreasing sequence $r_{n}$
with $r_{n}>f^{-1}(n),$ it holds $x\in \lim \sup T^{-n}(B(y,r_{n}))$.
\end{proposition}

\begin{proof}
The proof of this proposition is similar to the above proof of Theorem~\ref%
{mst}. Indeed we have that $\rho _{n}\leq f^{-1}(f(\rho_{n})) <
r_{\left\lfloor f(\rho_{n})\right\rfloor }$ and the proof follows as before.
\end{proof}

Interval Exchanges are particular bijective piecewise isometries which
preserve the Lesbegue measure. We refer to \cite{Bosher} for generalities on
this important class of maps. We only remark that Interval exchanges are not
hyperbolic and never mixing, hence Borel-Cantelli results about this class
of systems cannot come from speed of mixing arguments, as in \cite{T}. These
results will come from arithmetic arguments like in the rotation case. Let $%
T $ be some interval exchange. Let $\delta (n)$ be the minimum distance
between the discontinuity points of $T^{n}.$ We say that $T$ has the
property $\tilde{P}$ if it is ergodic and there is a\ constant $C$ and a
sequence $n_{k}$ such that $\delta (n_{k})\geq \frac{C}{n_{k}}.$

\begin{lemma}
\label{bosh9}(by \cite{Bosher}) The set of interval exchanges having the
property $\tilde{P}$ has full measure in the space of ergodic interval
exchange maps.
\end{lemma}

Now we can apply the above Proposition~\ref{aboveprop} to obtain the
following:

\begin{theorem}
If $T$ has the property $\tilde{P}$, (hence for typical interval exchanges)
there is a constant $K$ such that if $r_{n}$ is a decreasing sequence and $%
r_{n}\geq \frac{K}{n}$ eventually when $n$ is big enough, then the sequence $%
\{B(y,r_{n})\}$ has the BC property for almost each $y\in \lbrack 0,1]$.
\end{theorem}

\begin{proof}
In \cite{G} (in the proof of Theorem 9), it is proved that if $T$ has
property $\tilde{P}$ it holds that there is a sequence $\rho _{n}\rightarrow
0$ such that $\tau _{B(y,\rho _{n})}(x)\leq \frac{4}{C\rho _{n}}$ ($C$ is
the constant in the definition of property $\tilde{P}$ and may depend on $T$%
) for $x$ and $y$ ranging in positive measure sets $B$ and $B^{\prime }.$ We
now remark that, if $A_{n}$ is a decreasing sequence of sets, then $A=\lim
\sup T^{-n}(A_{n})$ is a forward invariant set, hence in an ergodic system
this set has either zero or full measure. By Proposition~\ref{aboveprop} we
have that if $y\in B^{\prime }$ and$\ r_{n}>\frac{4}{Cn}$ eventually, then $%
B\subset $ $\lim \sup T^{-n}(B(y,r_{n}))$. This implies $\mu (\lim \sup
T^{-n}(B(y,r_{n})))=1$ and that $\{B(y,r_{n})\}$ has the BC property.
Choosing $K>\frac{4}{C}$ we have the result for $y\in B^{\prime }.$

Let us consider a sequence $r_{n}$ such that $r_{n}\geq \frac{K}{n}$
eventually and $y$ such that $T(y)\in B^{\prime }$. Since $r_{n}$ is
decreasing this implie that the sequence $r_{n-1}$ is such that $r_{n-1}\geq 
\frac{K}{n}$ eventually (here we set $r_{-1}=1$) and then by what is proved
above $\{B(T(y),r_{n-1})\}$ is a BC sequence.

Now, if $\{B(T(y),r_{n-1})\}$ is a BC sequence of decreasing balls and both $%
y$ and $T(y)$ are not discontinuity points then also $\{B(y,r_{n})\}$ is a
BC sequence. This is true because $T^{-1}$ is an isometry from a small
neighborhood of $T(y)$ to a small neighborhood of $y$. This proves the
required result for each $y\in B^{\prime }\cup T^{-1}(B^{\prime })$ and the
result follows by the ergodicity of $T^{-1}$.
\end{proof}

Hence not only in rotations (by the results cited at the beginning of this
section), but also in a full measure set of interval exchanges we have that
a large class of decreasing sequences of centered balls have the BC property.


\begin{thebibliography}{99}
\bibitem{BGI} C. Bonanno, S. Isola and S. Galatolo, \emph{Recurrence and
algorithmic information}, Nonlinearity \textbf{17} (2004), no. 3, 1057--1074.

\bibitem{Bosher} M.~Boshernitzan, \emph{Quantitative recurrence results},
Invent. Math. \textbf{113} (1993), no.~3, 617--631.

\bibitem{bosh2} M. Boshernitzan, \emph{A condition for minimal interval
exchange maps to be uniquely ergodic}. Duke Math. J. \textbf{52} (1985),
no.~3, 723--752.

\bibitem{CG06} T. Carletti and S. Galatolo, \emph{Numerical Estimates of
dimension by recurrence and waiting time}, Physica A \textbf{364} (2006),
120--128.

\bibitem{CK} N.~Chernov and D.~Kleinbock, \emph{Dynamical {Borel-Cantelli}
lemmas for {Gibbs} measures}, Israel J. Math. \textbf{122} (2001), 1--27.

\bibitem{Dol} D.~Dolgopyat, \emph{Limit theorems for partially hyperbolic
systems}, Trans. Amer. Math. Soc. \textbf{356} (2004), 1637--1689.

\bibitem{L} N. Haydn, Y. Lacroix and S. Vaienti \emph{Hitting and return
times in ergodic dynamical systems}, Ann. Probab. \textbf{33} (2005), no.~5,
2043--2050.

\bibitem{HV} R. Hill, S. Velani \emph{Ergodic theory of shrinking targets, }%
Invent. Math. \textbf{119} (1995) 175--198.

\bibitem{kimkim} C. Kim and D.H. Kim, \emph{On the law of logarithm of the
recurrence time}, Discrete Contin. Dyn. Syst. \textbf{10} (2004), no.~3,
581--587.

\bibitem{kim2} D.H. Kim \emph{The dynamical Borel-Cantelli lemma for
interval maps}, Discrete Contin. Dyn. Syst. \textbf{17} (2007), no.~4,
891--900.

\bibitem{Kimnew} D.H. Kim \emph{The shrinking target property of irrational
rotations}, Nonlinearity \textbf{20} (2007), no.~7, 1637--1643.

\bibitem{Non2003} D.H. Kim and B.K. Seo, \emph{The waiting time for
irrational rotations}, Nonlinearity \textbf{16} (2003), no.~5, 1861--1868.

\bibitem{F} B. Fayad \emph{Two remarks on the shrinking target property}, J.
London Math. Soc. (to appear)

\bibitem{FMP} J.L. Fernandez, M.V Melian and D. Pestana, \emph{Quantitative
mixing results and inner functions}, Math. Ann. \textbf{337} (2007), no.~1,
233--251.

\bibitem{G05} S. Galatolo \emph{Dimension via waiting time and recurrence},
Math. Res. Lett. \textbf{12} (2005) no.~3, 377--386.

\bibitem{G} S. Galatolo \emph{Hitting time and dimension in Axiom A systems,
interval excanges and an application to Birkoff sums}, J. Stat. Phys. 
\textbf{123} (2006), 111--124

\bibitem{G07} S. Galatolo \emph{Dimension and hitting time in rapidly mixing
systems}, Math. Res. Lett. (in print)

\bibitem{Kurzweil} J.~Kurzweil, \emph{On the metric theory of inhomogeneous
diophantine approximatiions}, Studia Math. \textbf{15} (1955), 84--112.

\bibitem{K} Kleinbock D, \emph{Ergodic Theory on Homogeneous Spaces and
Metric Number Theory, }Encyclopedia of Complexity and Systems Science,
Springer Verlag, 2007

\bibitem{Ph} W.~Philipp, \emph{Some metrical theorems in number theory},
Pacific J. Math. \textbf{20} (1967), no.~1, 109--127.

\bibitem{Sh} P.~Shields, \emph{Waiting times: positive and negative results
on the Wyner-Ziv problem}, J. Theoret. Probab. \textbf{6} (1993), no. 3,
499--519.

\bibitem{T} D. Tasche, \emph{On the second Borel-Cantelli lemma for strongly
mixing sequences of events}, J. Appl. Probab. 34 (1997), no.~2, 381--394.
\end{thebibliography}
\end{document}